 
\baselineskip=14pt
\parskip=10pt

\def\epsilon{\varepsilon}

\def\halmos{\hbox{\vrule height0.15cm width0.01cm\vbox{\hrule height
 0.01cm width0.2cm \vskip0.15cm \hrule height 0.01cm width0.2cm}\vrule
 height0.15cm width 0.01cm}}
\font\eightrm=cmr8  
\font\eighttt=cmtt8
\magnification=\magstephalf

\parindent=0pt
\overfullrule=0in
\bf
\qquad Self Avoiding Walks, The Language\footnote{$^1$}
{\eightrm This paper is dedicated to the memory of my father
Yehudah Heinz Zeilberger(1915-1994), who spoke, read, and wrote fluently
in seven (natural) languages.}
of Science, and Fibonacci Numbers
 
\it
\quad \quad \quad
\qquad \qquad \qquad \qquad \qquad \qquad Doron ZEILBERGER \footnote{$^2$}
{\eightrm  \raggedright
Math Dept., Temple Univ.,
Philadelphia, PA 19122.
{\eighttt zeilberg@math.temple.edu .}
Supported in part by the NSF.
}
 
\rm
Science is a language.  In fact, science {\it is} language.
This was shown brilliantly by Xavier Viennot
and his {\it Ecole Bordelaise} (e.g. [V][DV][B]). Viennot, Maylis Delest,
and their disciples code {\it animals} and other {\it physical} creatures
in terms of {\it algebraic} (context-free) {\it languages}, by using
the so-called {\it Schutzenberger methodology}  (which Marco Schutzenberger
prefers to call the {\it DSV} (Dyck-Schutzenberger-Viennot) methodology.)
 
In this note, I  use this philosophy, or rather a juvenile
version of it, to enumerate {\it self avoiding walks} in the
(discrete) region $\{ 0,1 \} \times [ - \infty , \infty ]$, by encoding
these walks in terms of {\it words} in a certain {\it rational} (``finite-
automata'') language, that I call the  ``$UL^{*}IU'$ language'', 
and by describing its syntax.
 
A {\it self avoiding walk} (saw) in the two-dimensional (square) lattice
is a finite sequence of {\it distinct} lattice points 
$[(x_0 , y_0 )=(0,0) , (x_1 , y_1 ) , \dots , ( x_n , y_n )]$, such that
for all $i$, $(x_i , y_i )$ and $( x_{i+1} , y_{i+1} )$  are
{\it nearest neighbors}. The nearest neighbors of a point
$(a,b)$  are the four points $(a+1,b)$,$(a-1,b)$,$(a,b+1)$,$(a,b-1)$.
The problem of finding the exact, and even asymptotic, value of $a_n$,
the number of saws with $n$ steps, is wide open, and
presumably very difficult. The analogous problem in dimensions
$>4$, for the asymptotics, was recently solved brilliantly by
Hara and Slade[HS], and beautifully exposited in the masterpiece
by Madras and Slade[MS].
 
When a problem seems intractable, it is often a good idea to try
to study ``toy'' versions of it in the hope that
as the toys become increasingly larger and more sophisticated, 
they would metamorphose, in the limit, to the {\it real thing}. 
That was essentially
Lars Onsager's[O](see[T]) way of solving the Ising model.
Onsager first solved the ``finitary'' Ising problem
in a strip of finite-width, that turned out to be a
finite (definite) sum, that miraculously converged, a la Riemann-Integral,
to a certain definite integral.
 
Alm and Janson[AJ] had a similar idea of approaching general saws
by studying saws on strips $[-L,M] \times [ - \infty , \infty ]$,
with $L$ and $M$ finite. Saws,  when viewed ``dynamically'', are the
epitome of non-Markovianess. In [AJ] it was observed that when
saws are viewed ``statically'',  and restricted
to a strip, they can be described as Markov Chains.
A  saw can be viewed statically, since the path a self-avoiding drunkard 
makes uniquely determines her (or his) history. 
The general saw can be similarly viewed as a ``Markov chain'', but this
time the number of states is infinite. Since it is much easier
to describe the states then the saws themselves, there is some hope
that, by replacing the transition matrices by suitable
operators on some Hilbert (or whatever) space, this approach
will conquer the general problem. Only now
we transcend the rational, and even algebraic paradigms, into the
holonomic paradigm and beyond.
 
Alm and Jansen's[AJ] motivation was numerical. They wanted to find
lower bounds for the {\it connective constant}, $\mu$
($:= \lim_{n \rightarrow \infty} a_n^{1/n }$), by computing
the corresponding connective constants for saws in increasingly wider strips.
These turn out to be eigenvalues of matrices with integer entries,
and hence algebraic numbers.
 
Myself, I care little for {\it real}, floating-point numbers.
Being a discretian, I strive to get the {\it exact} answer.
The theorem below gives an exact enumeration of
$n$-step saws in the strip $\{ 0 , 1 \} \times [ - \infty , \infty ]$.
More interesting than the result is the 
{\it linguistic method of proof}, that would hopefully generalize.
 
{\bf Theorem}: The number, $a_n^{(2)}$, of $n-$step saws in the
strip $\{ 0,1 \} \times [ - \infty , \infty ]$
is given by $a_0^{(2)}=1$, $a_1^{(2)}=3$, and for $n > 1$, by
 
$$
a_n^{(2)} \, = \, 8 F_n \,- \,{{n} \over {2}} (1+ (-1)^n) \, -  \, 2(1- (-1)^n)
 \quad .
$$

{\bf Proof}: We assume that readers are familiar with the language
of {\it generatingfunctionlogy}[W].\footnote{$^3$}
{\eightrm Buy your own copy today! It would cost you less than
1/4 cent per day (Gian-Carlo Rota, in the ``Bulletin for
Mathematics books and software'', states that ``this book is good for the
next fifty years''.)}
From now on, let gf stand for ``(ordinary) generating function''.
 
Any saw in  $\{ 0,1 \} \times [ - \infty , \infty ]$
has the form $U L^{*} I U'$, where the meanings of $U, L, I, U'$
are as follows. (Steps in the right, left, up, and down direction
will be denoted by $r,l,u$, and $d$ respectively. For example
the walk $(0,0),(0,1),(0,2),(1,2),(1,1)$ will be coded as
$uurd$. Also $d^i$ means $d d \dots d$, where $d$ is repeated
$i$ times.\footnote{$^4$}{ \eightrm
Puzzle: What word in the English
language has the largest number of double letters?
Ans: $su b^2 o^2 k^2 e^2 per$})
 
(i) $U$ is a U-turn: $d^i r u^i$, with $i \geq 0$ ($i=0$
corresponds to a degenerate U-turn) (gf= $t/(1-t^2)$), or
nothing (gf=$1$). Total gf for this part is
$1+ t/(1-t^2)$.
 
(ii) $L^{*}$: Any number of (upside-down)$L$s (or $\Gamma$s), interlaced with 
upside-down-dyslectic $L$s. A single $L$ is either $u^i l$  or
$ u^{i} r$ ($i \geq 1$). The gf of a single $L$ is
$t^2 / (1- t)$ , and hence that of $L^{*}$ is
$1/(1- [t^2 / (1- t )] )= (1-t)/(1-t - t^2 )$.
(Philofibonaccist rejoice!)
 
(iii) an $I$, or nothing: $u^i , i \geq 0$. Its gf is $1/(1-t)$.
 
(iv) A final, modified U-turn, that I call $U'$, which is
$u^{i+1} l d^i$, or $u^{i+1} r d^i$, $i \geq 1$
(gf= $t^4/(1-t^2)$), or nothing (gf=$1$). The total gf is
$1+ t^4/(1-t^2)$.
 
The gf for the combined words $UL^{*}IU'$ is thus:
 
$$
[ 1+  {{t} \over {(1-t^2)}} ] \cdot 
[ {{(1-t)} \over {(1-t - t^2 )}} ] \cdot
[{{1} \over {(1-t)}}] \cdot [ 1+ {{t^4} \over {(1-t^2)}} ] =
{{(1+t-t^2)(1-t^2+t^4)} \over {(1-t^2)^2 (1-t- t^2 )}} \quad.
$$. 
 
But this is only half of the story: the northbound walks.
By symmetry, the gf of the other half, the southbound walks, which are the
x-axis mirror-reflection of the first half, is the same.
But two walks have been counted twice: the $0$-step empty walk
(gf=$1$), and the  $1$-step walk $[(0,0),(1,0)]=r$ (gf=$t$). So the
final gf is twice the gf above, take away $1+t$, namely
 
$$
2 {{(1+t-t^2)(1-t^2+t^4)} \over {(1-t^2)^2 (1-t- t^2 )}} \, - \,
(1+t) \, = \,
{{1+2t - t^3 -t^4 + t^7 } \over {(1-t)^2 (1+t)^2 (1-t-t^2)} }
\quad.
$$
 
A partial-fraction decomposition (that Maple$^{TM}$ kindly performed
for me), followed by a Maclaurin expansion, yields the formula
for $a_n^{(2)}$. \halmos
 
{\bf A Shorter, more elegant, Semi-Rigorous, late-21st Century-Style Proof:}
Compute $a_n^{(2)}$ by direct enumeration for $0 \leq n \leq 15$, and then use
Salvy and Zimmerman's[SZ]  Maple package {\it gfun} to conjecture
the gf. Since we know a priori that this is a rational function,
that must be it. \halmos
 
To make this argument completely rigorous,
you would have to derive a priori bounds for the degrees of the numerator
and denominator of the gf, but {\it who cares ?}
 
The only possible advantage of the first proof is that it might
generalize to obtain the gfs, $\phi_r (t)$ for the
number of saws  in the strip $[-r,r] \times [ - \infty , \infty ]$,
for $r=1,2, \dots$. Of course, the expressions themselves will very soon
become unwieldy. More exciting is the prospect that one might
be able to find some kind of functional equation that
expresses $\phi_r (t)$ in terms of $\phi_{r-1} (t)$, or more
refined quantities, from which the divine quantity
$\phi (t) :=  \lim_{r \rightarrow \infty} \phi_r (t)$
could be looked at in the eyes, without being blinded.
{\it Amen ken yehi ratson}.
\medskip
{\bf REFERENCES}
 
[AJ] S.E. Alm and S. Janson, {\it Random self-avoiding walks
on one-dimensional lattices}, Comun. Statist.-Stochastic Models
{\bf 6}(1990), 169-212.
 
[B] M. Bousquet-M\'elou, ``{\it q-\'Enum\'eration de polyominos convexes}'',
Publications de LACIM, UQAM, Montr\'eal, 1991.
 
[DV] M.P. Delest and X.G. Viennot, {\it Algebraic languages and polyominoes
enumeration}, Theor. Comp. Sci. {\bf 34}(1984), 169-206.
 
[HS] T. Hara, and G. Slade, {\it The lace expansion for self-avoiding
walk in five or more dimensions}, Reviews in Math. Phys. {\bf 4}(1992),
235-327.
 
[MS] N. Madras and G. Slade, ``{\it The Self Avoiding Walk}'',
Birkhauser, Boston, 1993.
 
[O] L. Onsager, {\it Crystal Statistics, I. A two-dimensional model
with an order-disorder transition}, Phys. Rev. {\bf 65}(1944), 117.
 
[SZ] B. Salvy and P. Zimmerman, {\it gfun : a Maple package for the
manipulation of generating and holonomic functions in one variable},
to appear in ACM Transactions on Mathematical Software.
 
[T] C.J. Thompson, ``{\it Mathematical Statistical Mechanics}'',
Princeton University Press, Princeton, 1972.
 
[V] X.G. Viennot, {\it Probl\`emes combinatoire pos\'es par la 
physique statistique, S\'eminaire Bourbaki $n^{o}$} 626,
Asterisque {\bf 121-122}(1985), 225-246.
 
[W] H. Wilf, ``{\it Generatingfunctionology}''($2^{nd}$ edition), 
Academic Press, San Diego, 1994.
 
Yud beIyar  Heh' Tashnad (10 Iyar, 5754)
\bye